\documentclass[journal]{IEEEtran}

\usepackage{cite}
\usepackage{todonotes}
\usepackage{graphicx}
\usepackage{subfigure}
\usepackage{epstopdf}
\usepackage{url}
\usepackage{stfloats}
\usepackage{latexsym}
\usepackage{amsfonts}
\usepackage{amsmath}
\usepackage{amssymb}
\usepackage{array}
\usepackage{textcomp}
\usepackage{pxfonts}
\usepackage{txfonts}
\usepackage{cases}
\usepackage{psfrag}
\usepackage{epsfig}
\usepackage{setspace}
\usepackage{enumitem}
\usepackage{xcolor}
\usepackage{multirow}
\usepackage{multicol}
\usepackage{lmodern}

\usepackage[T1]{fontenc}
\usepackage{color}
\usepackage{float}
\usepackage{esint}

\makeatletter

\floatstyle{ruled}
\newfloat{algorithm}{tbp}{loa}
\providecommand{\algorithmname}{Algorithm}
\floatname{algorithm}{\protect\algorithmname}

\usepackage{fullpage}
\usepackage{algorithmic,algorithm}

\makeatother

\newtheorem{Definition}{Definition}
\newtheorem{Remark}{Remark}
\newtheorem{Lemma}{Lemma}

\newtheorem{Proposition}{Proposition}

\usepackage{todonotes}   

\newcommand{\R}{\ensuremath{\mathbb{R}}}
\newcommand{\N}{\ensuremath{\mathbb{N}}}

\newcommand{\Z}{\ensuremath{\mathbb{Z}}}

\newcommand{\cR}{\mathcal{R}}
\newcommand{\cL}{\mathcal{L}}

\newlist{aims}{enumerate}{1}
\setlist[aims,1]{
  label={Additive Faults, set~\arabic*:},
  leftmargin=*,
  align=left,
  labelsep=10mm,
}

\newlist{steps}{enumerate}{1}
\setlist[steps,1]{
  label={Step~\arabic*:},
  leftmargin=*,
  align=left,
  labelsep=10mm,
}

\begin{document}

\title{\Large Structural Identifiability Analysis of Fractional Order Models with Applications in Battery Systems}
\author{S.M.Mahdi Alavi, Adam Mahdi, Pierre E. Jacob, Stephen J. Payne, and David A. Howey
\thanks{S.M.M. Alavi was with the Energy and Power Group, Department of Engineering Science, University of Oxford. He is now with the Brain Stimulation Engineering Laboratory, Department of Psychiatry \& Behavioral Sciences, Duke University, Durham, NC 27710, USA. Email: mahdi.alavi@duke.edu}
\thanks{A. Mahdi and S.J. Payne are with the Institute of Biomedical Engineering, Department of Engineering Science, University of Oxford, Old Road Campus Research Building, Oxford, OX3 7DQ, United Kingdom. Emails:  \{adam.mahdi, stephen.payne\}@eng.ox.ac.uk}
\thanks{P.E. Jacob is with the Department of Statistics, Harvard University, Science Center 7th floor, 1 Oxford Street, Cambridge, MA 02138-2901, USA.  Email: pierre.jacob.work@gmail.com}
\thanks{D.A. Howey is with the Energy and Power Group, Department of Engineering Science, University of Oxford, Parks Road, Oxford, OX1 3PJ, United Kingdom. Email: david.howey@eng.ox.ac.uk}
}

\markboth{ }%
{Shell \MakeLowercase{\textit{et al.}}: Bare Demo of IEEEtran.cls
for Journals}
\maketitle

\begin{abstract}
This paper presents a method for structural identifiability analysis of fractional order systems by using the coefficient mapping concept to determine whether the model parameters can uniquely be identified from input-output data. The proposed method is applicable to general non-commensurate fractional order models. Examples are chosen from battery fractional order equivalent circuit models (FO-ECMs). 
The battery FO-ECM consists of a series of parallel resistors and constant phase elements (CPEs) with fractional derivatives appearing in the CPEs. The FO-ECM is non-commensurate if more than one CPE is considered in the model. Currently, estimation of battery FO-ECMs is performed mainly by fitting in the frequency domain, requiring costly electrochemical impedance spectroscopy equipment. This paper aims to analyse the structural identifiability of battery FO-ECMs directly in the time domain. It is shown that FO-ECMs with finite numbers of CPEs are structurally identifiable. In particular, the FO-ECM with a single CPE is structurally globally identifiable. 
\end{abstract}

\begin{IEEEkeywords}
System identification, Identifiability, Fractional order systems, Batteries.
\end{IEEEkeywords}

\IEEEpeerreviewmaketitle

\section{Introduction}
\label{sec:introduction}
Parameter estimation is very important for studying and understanding of systems \cite{Astroem1971, Young1981, Ljung2010,Ljung1987, Soderstrom1989, Zhu2001}. 
However, because of the specific model structure or inadequate data it might not be possible to infer some (or all) of the model parameters. 
\emph{Structural identifiability} analysis is a data-free concept, which determines whether the model parameters can be uniquely estimated from rich input-output data \cite{Bellman1970}.
\emph{Practical identifiability}, on the other hand, takes into account the practical aspects of the problem that come with real data including noise, bias and signal quality such as  the shape of excitation, its magnitude, length, frequency ranges, etc, \cite{Raue2009}. Thus, structural identifiability analysis is a prerequisite for practical identifiability.

Since the 1970's, several analytical techniques have been proposed for studying structural identifiability analysis based on Taylor series expansion \cite{Pohjanpalo1978, Chappell1990}, similarity transformations \cite{Vajda1989, Anstett2008, Meshkat2014, Mahdi2014, Glover1974, DistefanoIII1977, VanDenHof1998}, differential algebra \cite{Ljung1994,Audoly2001}, and Laplace transforms (transfer functions) \cite{Cobelli1980, Bellman1970, Nazarian2010}. However, almost all of the proposed techniques deal with ordinary differential equations with integer orders. To the best of our knowledge, only reference \cite{Nazarian2010} studies the identifiability of single-input single-output (SISO) fractional commensurate-order systems in the frequency domain. In the fractional commensurate-order systems, all the orders of derivation are integer multiples of a base order, \cite{Oustaloup1995, Chen2009, Caponetto2010, Monje2010}. It was shown that SISO fractional commensurate-order systems are poorly identifiable for small values of the base order, \cite{Nazarian2010}.

This paper presents an alternative methodology to determine the structural identifiability of SISO fractional-order (FO) systems. The method is applicable to both commensurate and non-commensurate models based in time domain. This is an advantage since typically, in real world applications, the data is given in the time domain.

Recently, there has been a significant interest in identifiability analysis of battery models. References \cite{Schmidt2010, Forman2012, Moura2014, DAmato2012} study the identifiability of electrochemical models and show that some parameters are not practically identifiable. In \cite{Sitterly2011, Rausch2013, Sharma2014} and \cite{Rothenberger2014}, the identifiability of conventional equivalent circuit models (ECMs) is addressed. Conventional ECMs only include resistors and standard capacitors based on the Randles circuit \cite{Randles1947}. Their dynamics can be expressed by ordinary differential equations with integer orders. A comprehensive survey of the conventional ECMs with integer orders has been given in \cite{Hu2012}. However, these may not accurately reflect the dynamic behaviour of real battery systems. In order to address these issues, constant phase elements (CPEs) with fractional-order dynamics are incorporated into ECMs, which results in fractional-order ECMs (FO-ECMs), \cite{Cole1941}. Compared to conventional ECMs, FO-ECMs represent distributed electrode processes more accurately, \cite{Barsoukov2005} and may give more insight into battery performance, which could be useful for monitoring and diagnostic purposes \cite{Troeltzsch2006, Richardson2014}.

The usefulness of a FO-ECM highly depends on the ability to estimate its parameters. Currently, this is done by fitting the model to frequency domain impedance spectra that are obtained through Electrochemical Impedance Spectroscopy (EIS), \cite{Troeltzsch2006, Macdonald1982, Boukamp1986}. However, parameter estimation of FO-ECMs directly from time-domain data is very appealing since conversion to the frequency domain may introduce bias in the estimation \cite{Alavi2015}.

In this paper, a structural identifiability analysis method based on the concept of coefficient map is employed, which is applicable to general non-commensurate fractional order models. The method is applied to study structural identifiability of battery FO-ECMs, which are non-commensurate if more than one CPE is considered in the model. It is shown that the structural identifiability of battery FO-ECMs depends on the solution of a set of nonlinear coupled equations. The number of equations to be solved equals the number of CPEs. It is shown that FO-ECMs with finite numbers of CPEs are structurally identifiable and the FO-ECM with a single CPE is structurally globally identifiable. 



\section{Model Structure}\label{sec:model}
In this section discrete-time state-space and transfer function models of fractional-order systems are derived.

A state-space representation of a SISO FO system is given by
\begin{eqnarray}\label{SS:c}
\begin{aligned}
\frac{d^{\alpha} x(t)}{dt^{\alpha}}&= \bar{A}(\beta)\, x(t) + \bar{B}(\beta)\,u(t)\\
 y(t) &= M(\beta)\, x(t)+ D (\beta)\,u(t)
\end{aligned}
\end{eqnarray}
where $x(t) \in \R^n$ is the state vector; $u(t) \in \R$ and $y(t) \in \R$ are input and output signals, respectively; $\bar{A}(\beta)\in \R^{n \times n}$, $\bar{B}(\beta) \in \R^{n \times 1}$, $M(\beta) \in \R^{1 \times n}$ and $D(\beta) \in \R$ are system matrices which depend on the parameter vector $\beta$ to be identified. Moreover,
\begin{equation}\label{fracvec}
\frac{d^{\alpha} x(t)}{dt^{\alpha}} = \Big[\frac{d^{\alpha_1} x_1(t)}{dt^{\alpha_1}},\ldots,\frac{d^{\alpha_n} x_n(t)}{dt^{\alpha_n}}\Big]^\top
\end{equation}
is the vector of fractional-order derivatives with the unknown fractional-orders  $\alpha_i\in(0,1)$, $i=1,\ldots,n$.

\begin{Definition}
\label{Def:commensurateFOsys}A FO system is said to be commensurate if $i\in\{1,\ldots n\},~\exists \rho \in \N,~\mbox{such~that}~\alpha_i =  \rho \alpha$, where $\alpha \in \R$; otherwise it is said to be non-commensurate, \cite{Oustaloup1995, Chen2009, Caponetto2010, Monje2010}. \hfill{$\square$}
\end{Definition}

This paper considers general SISO non-commensurate FO systems.



A discrete-time representation of the fractional differentiation operator $d^{\alpha} x(t)/dt^{\alpha}$ is often given by the Gr\"{u}nwald-Letnikov approximation \cite{Oustaloup1995}
\begin{align}
\nonumber
&\mbox{diag}\{T_s^{\alpha_1},\ldots,T_s^{\alpha_n}\} \frac{d^{\alpha} x(kT_s)}{dt^{\alpha}}=\\ 
\label{G-L derivative} &
\displaystyle \sum_{j=0}^{k+1}(-1)^j\mbox{diag}\left\{\binom{\alpha_1}{j},\cdots,\binom{\alpha_n}{j}\right\} x((k+1-j)T_s),
\end{align}
where $T_s$ is the sample time, $k\in \Z^+$ is the time index,  $\mbox{diag}\{\cdot\}$ denotes the diagonal matrix and $\binom{\alpha_i}{j}$ is the binomial coefficient given by
\begin{equation}\label{Gamma_func}
\binom{\alpha_i}{j}=\frac{\Gamma(\alpha_i+1)}{\Gamma(j+1)\Gamma(\alpha_i+1-j)},
\end{equation}
where, $\Gamma(\cdot)$ denotes the gamma function
\[
\Gamma(\alpha_i)=\int_{0}^{\infty} z^{\alpha_i-1} e^{-z}dz,~ \text{for} ~ \alpha_i\in\mathbb{C}~ \text{with}~  \Re(\alpha_i)>0.
\]
For the sake of simplicity, $T_s$ is omitted from the argument and $x(k+1-j)$ is written as $x_{k+1-j}$ hereafter. The substitution of the Gr\"{u}nwald-Letnikov approximation into equation \eqref{SS:c} gives
\begin{align}
\nonumber
& \mbox{diag}\{T_s^{-\alpha_1},\ldots,T_s^{-\alpha_n}\} \times\\
\nonumber & \displaystyle  \sum_{j=0}^{k+1}\left((-1)^j\mbox{diag}\left\{\binom{\alpha_1}{j},\cdots,\binom{\alpha_n}{j}\right\}x_{k+1-j} \right)\\
\label{aux_ss1} & \hspace{7em}
= \bar{A}(\beta) x_k +\bar{B}(\beta) u_k.
\end{align}
Multiplying \eqref{aux_ss1} by $\mbox{diag}\{T_s^{\alpha_1},\ldots,T_s^{\alpha_n}\}$ and extracting  $x_{k+1}$ from the summation gives
\begin{align}
\nonumber & x_{k+1}=\\
\nonumber & \left(\mbox{diag}\{\alpha_1,\cdots,\alpha_n\}+\mbox{diag}\{T_s^{\alpha_1},\cdots,T_s^{\alpha_n}\}\bar{A}(\beta)\right) x_k+\\
\nonumber & \mbox{diag}\{T_s^{\alpha_1},\cdots,T_s^{\alpha_n}\}\bar{B}(\beta) u_k -\\
\label{aus_ss2}& \displaystyle \sum_{j=2}^{k+1}\left((-1)^j\mbox{diag}\left\{\binom{\alpha_1}{j},\cdots,\binom{\alpha_n}{j}\right\}x_{k+1-j} \right).
\end{align}

Thus a discrete-time state-space model, in compact form, can be written as
\begin{eqnarray}\label{SS:d}
\begin{aligned}
 x_{k+1}&=\displaystyle \sum_{j=0}^{k} A_j(\beta,\alpha)\, x_{k-j}+B(\beta,\alpha) u_k \\
\label{output equation discrete}
y_k&=M(\beta) x_k + D(\beta) u_k.
\end{aligned}
\end{eqnarray}
with
\begin{eqnarray}\label{SS:d-matrices}
\begin{aligned}
&\alpha=\left[\begin{array}{ccc} \alpha_1 ~ \cdots ~ \alpha_n \end{array}\right]\\
& A_0=\mbox{diag}\{\alpha_1,\cdots,\alpha_n\}+\\
& \hspace{3em}\mbox{diag}\{T_s^{\alpha_1},\cdots,T_s^{\alpha_n}\}\bar{A}(\beta)\\
& A_j=-(-1)^{j+1}\mbox{diag}\left\{\binom{\alpha_1}{j+1},\cdots,\binom{\alpha_n}{j+1}\right\},\\
&\hspace{3em} \mbox{~for~} 1 \leq j\\
& B =\mbox{diag}\{T_s^{\alpha_1},\cdots,T_s^{\alpha_n}\}\bar{B}(\beta).
\end{aligned}
\end{eqnarray}

Finally, the transfer function model structure parametrised by the unknown parameters
\begin{align}
\theta=[\beta ~ \alpha]
\end{align} is given by
\begin{align}
\nonumber
H(z,\theta)=&M(\theta)\left(zI-\displaystyle \sum_{j=0}^{T}z^{-j}A_j(\theta)\right)^{-1}B(\theta)+\\ \label{TF-MS-theta} & \hspace{2em} D(\theta),
\end{align}
where $I$ is the $n \times n$ identity matrix and $T$ represents the data length, i.e., the number of samples.

\begin{Remark}
The state-space model \eqref{SS:d} implies that $x_{k+1}$ depends on all the past states, $x_0$ up to $x_k$. This means that FO systems are  non-Markov\footnote{In Markov system, $x_{k+1}$ can be written as functions of $x_k$ and inputs.}. From the transfer function perspective this means that the order of a FO system's transfer function equals the data length (equation \eqref{TF-MS-theta}). These are the main distinguishing features of FO systems which make their analysis and identification challenging. \hfill{$\square$}
\end{Remark}



\section{Structural Identifiability}\label{sec:strident}
The goal of this section is to introduce rigorously the concept of structural identifiability, which will be applied in the following section to FO-ECMs.


\begin{Definition}\label{def:ide}(\cite{Ljung1987})
Consider a model $\mathcal{M}$ with the  transfer function $H(z,\theta)$, parametrised by $\theta$, where $\theta$ belongs to an open subset $\mathcal{D}_{\theta_{\mathcal{M}}} \subset \R^q$, and  consider the equation
\begin{equation}\label{TFd:ide}
H(z,\theta)=H(z,\theta^\ast),\qquad \text{for almost all } z,
\end{equation}
where $\theta, \theta^\ast\in\mathcal{D}_{\theta_{\mathcal{M}}}$. Then,  model $\mathcal{M}$ is said to be
\begin{itemize}
\item[-] \emph{globally identifiable} if \eqref{TFd:ide} has a unique solution in $\mathcal{D}_{\theta_{\mathcal{M}}}$,
\item[-] \emph{identifiable} if \eqref{TFd:ide} has a finite number of solutions in $\mathcal{D}_{\theta_{\mathcal{M}}}$,
\item[-] \emph{unidentifiabile} if \eqref{TFd:ide} has an infinite number of solutions in $\mathcal{D}_{\theta_{\mathcal{M}}}$.
\end{itemize}
\hfill{$\square$}
\end{Definition}

Instead of using the above definition of structural identifiability, which appears to be a standard definition in the engineering literature, it might be more convenient to use the concept of coefficient map, which will now be reviewed. 

\begin{Definition}
Consider the following monic transfer function\footnote{The coefficient of the highest order term in the denominator is 1.}
\begin{equation}\label{general-tfms}
H(z,\theta)=\frac{f_{n_f}(\theta)z^{n_f}+f_{n_f-1}(\theta)z^{n_f-1}+\cdots+f_{0}(\theta)}{z^{n_g}+g_{n_g-1}(\theta)z^{n_g-1}+\cdots+g_{0}(\theta)}
\end{equation}
where $\theta=\big[\begin{array}{ccc} \theta_1 & \cdots & \theta_q \end{array}\big]\in\mathcal{D}_{\theta_{\mathcal{M}}}$ is the parameter vector, and $\mathcal{D}_{\theta_{\mathcal{M}}}\subset \R^q$ is an open set. The \emph{coefficient map} $\mathcal{C}_{\mathcal{M}}:\mathcal{D}_{\theta_{\mathcal{M}}} \rightarrow \R^{n_f+n_g+1}$ is defined as
\begin{equation}\label{cm}
\mathcal{C}_{\mathcal{M}}(\theta)=(f_{n_f}(\theta),\cdots,f_{0}(\theta),g_{n_g-1}(\theta),\cdots, g_{0}(\theta)).
\end{equation}
\hfill{$\square$}
\end{Definition}

The following lemma illustrates the applicability of the coefficient map for studying the structural identifiability.

\begin{Lemma}[\cite{Mahdi2014, Meshkat2014}]\label{lem:Coefmap}
Consider a model $\mathcal{M}$ with the transfer function \eqref{general-tfms} and the associated coefficient map \eqref{cm}. Then model $\mathcal{M}$ is
\begin{itemize}
\item[-]  \emph{globally identifiable} if the coefficient map $\mathcal{C}_{\mathcal{M}}$ is one-to-one,
\item[-] \emph{identifiable} if the coefficient map $\mathcal{C}_{\mathcal{M}}$ is many-to-one,
\item[-] \emph{unidentifiable} if the coefficient map $\mathcal{C}_{\mathcal{M}}$ is infinitely many-to-one.
\end{itemize}
\hfill{$\square$}
\end{Lemma}

An important special case of \eqref{general-tfms} is when the coefficients of the transfer function form the parameter vector, in which case the identifiability is given by the following lemma.
\begin{Lemma}\label{lem:cop}
Consider the following transfer function
\begin{align}
\label{Coprime-lemma-TF}
H(z,\theta)=\frac{f_{n_f}z^{n_f}+f_{n_f-1}z^{n_f-1}+\cdots+f_{0}}{z^{n_g}+g_{n_g-1}z^{n_g-1}+\cdots+g_{0}},
\end{align}
where the parameter vector only consists of the coefficients of the numerator and denominator, i.e.  $\theta=[ f_{n_f}, \cdots,f_0, g_{n_{g-1}},\cdots, g_0]^\top$. Then, the model structure \eqref{Coprime-lemma-TF} is globally identifiable, see Section 4.6 of \cite{Ljung1987}. \hfill{$\square$}
\end{Lemma}

The above lemma is only applicable if the parameter vector is formed by the coefficients of the transfer function. If the coefficients of the transfer function are themselves functions of some parameter vector $\theta$, then the concept of reparametrisation can be used to study the identifiability of the model.


\begin{Definition}\label{def:Rep}
Consider a model structure $\mathcal{M}$ with the parameter vector $\theta_\mathcal{M}$ belonging to the open subset $\mathcal{D}_{\theta_\mathcal{M}} \subset \R^{q}$. A \emph{reparametrisation} of the model structure $\mathcal{M}$  with the coefficient map $\mathcal{C}_\mathcal{M}$ is a map $\mathcal{R}: \mathcal{D}_{\theta_\mathcal{M}} \to \R^{p}$ such that
\begin{equation}
Im\, (\mathcal{C}_{\mathcal{M}} \circ \mathcal{R}) = Im\, (\mathcal{C}_{\mathcal{M}}),
\end{equation}
where $Im$ denotes the image of the map, and `$\circ$' denotes composition.\hfill{$\square$}
\end{Definition}

\medskip
Based on the above definitions and lemmas, the process of determining the structural identifiability of a model $\mathcal{M}$ can be divided into three general steps:
\begin{algorithmic}[1]
	\STATE Compute the transfer function of model $\mathcal{M}$.
	\STATE Determine the corresponding coefficient map $\mathcal{C}_{\mathcal{M}}$.
	\STATE If $\mathcal{C}_{\mathcal{M}}$ is one-to-one then model $\mathcal{M}$ is globally identifiable; if $\mathcal{C}_{\mathcal{M}}$ is finitely many-to-one then $\mathcal{M}$ it is identifiable; finally if $\mathcal{C}_{\mathcal{M}}$ is infinitely many-to-one then model  $\mathcal{M}$ is unidentifiable.
\end{algorithmic}

\section{Structural Identifiability in Battery Systems}\label{sec:batsi}
In this section, the proposed identifiability analysis method is applied to battery models and results are discussed. First, a simple integer-order derivative example is considered in order to facilitate the reading of the more involved fractional-order case. 

\subsection{Integer-order models}

\label{exRRCstd}
\begin{figure}
\centering
\includegraphics[scale=0.8]{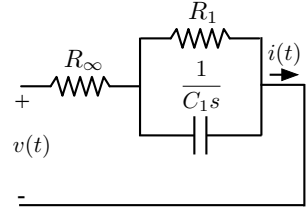}
\caption{The Randles circuit with standard capacitor.}
\label{Fig:StdRRC}
\end{figure}

Consider the electric circuit as shown in Figure\,\ref{Fig:StdRRC}. This circuit was proposed by Randles in 1947 \cite{Randles1947} for modeling the kinetics of rapid electrode reactions. Since then, the model has become the basis for studying  various electrochemical energy storage systems such as batteries, fuel cells and supercapacitors, \cite{Alavi2015a}. The resistor $R_\infty$ in this circuit models the battery ohmic resistance. The resistor $R_1$ and the capacitor $C_1$ denote diffusion processes or the battery charge transfer resistance (CTR) and double layer (DLC) capacitance, respectively. 

\begin{Proposition}
{\it The Randles model given in Figure\,\ref{Fig:StdRRC} with the parameter vector $\theta=[R_\infty, R_1, C_1]$ is structurally globally identifiable. }
\end{Proposition}

\smallskip

The above proposition will be shown in three steps.

\smallskip

{\em Step 1:} By Kirchhoff's laws, it is simple to show that a transfer function of the circuit, parameterised by $\theta$ is
\begin{align}
\label{exRRCstd:tf}
H(z,\theta)=\frac{f_1(\theta)z+f_0(\theta)}{z+g_0(\theta)},
\end{align}
where
\begin{eqnarray}\label{exRRCstd:tfcoef}
\begin{aligned}
& f_1(\theta)=R_\infty\\
& f_0(\theta)=-R_\infty(1-\frac{T_s}{R_1C_1})+\frac{T_s}{C_1}\\
& g_0(\theta)=-(1-\frac{T_s}{R_1C_1}).
\end{aligned}
\end{eqnarray}

{\em Step 2:} The coefficient map associated with the model is given by
\[
\mathcal{C}: \theta \to \big(f_1(\theta),f_0(\theta),g_0(\theta)\big).
\]

\smallskip

{\em Step 3:} For global identifiability it is sufficient to show that the coefficient map is one-to-one. Since the above function is invertible with the inverse given by
\[
\mathcal{C}^{-1}: (f_1,f_0,g_0) \to \Big(f_1, \frac{T_s}{C_1(1+g_0)}, \frac{T_s}{f_0-f_1g_0}\Big),
\]
the model \eqref{exRRCstd:tf} is globally identifiable. \hfill{$\square$}

\smallskip

More details about the structural and practical identifiability of the Randles circuit and its generalised topology are given in \cite{Alavi2015a}.


\begin{figure}
\centering
\includegraphics[scale=.8]{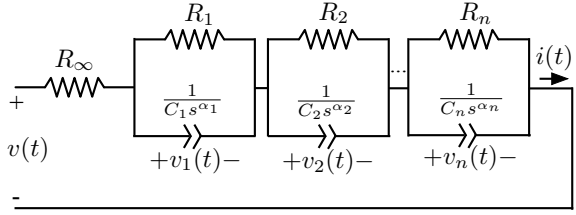}
\caption{The general battery electrochemical impedance spectroscopy model.}
\label{Fig:EISgeneral}
\end{figure}

\subsection{Fractional-order models}
In this section, we study the identifiability of FO-ECMs directly from time-domain data. A general impedance schematic of the EIS FO-ECM is shown in Figure \ref{Fig:EISgeneral}. There are two main differences between the EIS FO-ECM Figure \ref{Fig:EISgeneral} and the Randles circuit Figure \ref{Fig:StdRRC}. In the EIS FO-ECM, more than one parallel pair is seen. Each parallel pair is employed to model the battery  processes over a certain frequency range. The Randles model can also be developed by adding more parallel pairs as discussed in \cite{Alavi2015a}. The number of parallel pairs depends on the required accuracy for the frequency domain fitting of impedance spectra.

The second and main difference is that the impedance of the capacitor in the Randles circuit is given by the integer derivative, while in the FO-ECMs fractional derivatives are applied. As mentioned in Section \ref{sec:introduction}, these elements are referred to as constant phase elements (CPEs), \cite{Cole1941}. CPEs model diffusion processes (or CTRs/DLCs) more accurately as shown in \cite{Alavi2015}. The impedance of the i-th CPE is given by:
\begin{align}
\label{EIS-CPE}
Z_{CPE_i}(s)=\frac{1}{C_is^{\alpha_i}},
\end{align}
where $C_i$ is a constant, $s$ is the Laplace operator and $\alpha_i$ ($0 > \alpha_i > 1$) is the exponent value. The dimension of $C_i$ is $\mbox{Fcm}^{-2}s^{\alpha_i-1}$ \cite{Jorcin2006}. In low frequency ranges the impedance frequency response may show constant phase behaviour such that the associated parallel resistor can be considered as an open circuit. This is referred to as Warburg term in the literature \cite{Barsoukov2005}.

It should be noted that there are techniques that approximate CPEs with ideal capacitors \cite{Plett2004}, or by a series connection of numerous R-C pairs, \cite{Andre2011, Birkl2013, Hu2011}.

The parameter vector associated with Figure \ref{Fig:EISgeneral} is defined as follows:
\begin{eqnarray}\label{gen-ecm-parvec}
\begin{aligned}
\theta=[ R_\infty ,R_1 ,\cdots , R_n , C_1 , \cdots , C_n , \alpha_1 ,\cdots , \alpha_n].
\end{aligned}
\end{eqnarray}
By defining the voltage across the CPEs as the state variables,
\begin{align}
x \triangleq \left[\begin{array}{ccc} v_1 & \cdots & v_n\end{array}\right]^\top,
\end{align}
and by Kirchhoff's laws, it is easy to show that $A_j$, $B$, $M$ and $D$ in the state-space model \eqref{SS:d} are given by:

\begin{eqnarray}\label{gen-ecm-ssd}
\begin{aligned}
& A_j(\theta)=\mbox{diag}\left\{a_{1,j}(\theta), \cdots, a_{n,j}(\theta)\right\}\\
& B(\theta)=\left[\begin{array}{ccc} b_1(\theta) & \cdots & b_n(\theta)\end{array}\right]^\top\\
& M(\theta)=\left[\begin{array}{ccc} m_1 & \cdots & m_n\end{array}\right]\\
& D(\theta)=d(\theta),
\end{aligned}
\end{eqnarray}
with
\begin{eqnarray}\label{gen-ecm-params}
\begin{aligned}
& a_{i,0}(\theta)= \alpha_i-\frac{T_s^{\alpha_i}}{R_iC_i}\\
& a_{i,j}(\theta)= -(-1)^{j+1}\binom{\alpha_i}{j+1}\\
& b_i(\theta)=\frac{T_s^{\alpha_i}}{C_i}\\
& m_i=1\\
& d(\theta)=R_\infty\\
& \text{for~}i=1,\cdots,n \mbox{~and~} j=1,2,\cdots,T,
\end{aligned}
\end{eqnarray}
where $T$ is the data length. 

By using \eqref{TF-MS-theta}, a transfer function model structure of the circuit FO-ECM of Figure \ref{Fig:EISgeneral} is given by:
\begin{align}
\label{TF-gecm}
H(z,\theta)=d(\theta)+\displaystyle \sum_{i=1}^{n}\frac{m_ib_i(\theta)z^{T}}{z^{T+1}- \sum_{j=0}^{T} a_{i,j}(\theta)z^{T-j}}.
\end{align}

The structural identifiability of FO-ECMs with respectively one and two parallel R-CPE pairs is studied in the following sub-sections. The results provide useful insight into the general case with $n$ CPEs. From a practical perspective, the examples under consideration in the rest of the section are sufficient to represent the basic dynamic behaviour of a large number of electrochemical systems, \cite{Alavi2015, Waag2013}.

\subsubsection{$R_\infty - R_1||\frac{1}{C_1 s^{\alpha_1}}$ circuit}\label{exRRC}

\begin{figure}[ht]
\centering
\includegraphics[scale=0.7]{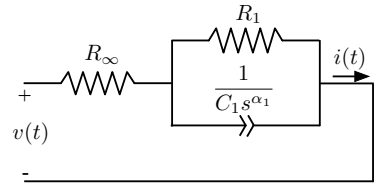}
\caption{The $R_\infty - R_1||\frac{1}{C_1 s^{\alpha_1}}$ circuit.}
\label{Fig:fracR-RC}
\end{figure}

Consider the EIS FO-ECM with a single parallel R-CPE in series with an ohmic resistor as shown in Figure \ref{Fig:fracR-RC}.  

\medskip

\begin{Proposition}
{\it The FO-ECM shown in Figure \ref{Fig:fracR-RC} with the parameter vector $\theta=[R_\infty,R_1, C_1, \alpha_1]$  is structurally globally identifiable.}
\end{Proposition}

\smallskip

The proposition will be shown in three steps.

\smallskip

{\em Step 1:} Using \eqref{TF-gecm}  the transfer function is found to be
\begin{align}\label{exRRC:tf1}
H(z,\theta)=d(\theta)+\frac{b_1(\theta)z^{T}}{z^{T+1}- \sum_{j=0}^{T} a_{1,j}(\theta)z^{T-j}},
\end{align}
where the relationships between the transfer function coefficients and the model parameters are given in \eqref{gen-ecm-params}. In order to compute the coefficient map,  \eqref{exRRC:tf1} is written as a monic rational function
\begin{align}
\label{exRRC:tf2}
H(z,\theta)=\frac{f_{T+1}(\theta)z^{T+1}+\cdots+ f_{0}(\theta)}{z^{T+1}+g_{T}(\theta)z^{T}+\cdots+ g_{0}(\theta)},
\end{align}
where
\begin{eqnarray}\label{exRRC:tfcoef}
\begin{aligned}
& f_{T+1}(\theta)=d(\theta)\\
& f_{T}(\theta)=b_1(\theta)-a_{1,0}(\theta)d(\theta)\\
& f_{T-j}(\theta)=-a_{1,j}(\theta)d(\theta), \mbox{~for~}1\leq j \leq T\\
& g_{T-j}(\theta)=-a_{1,j}(\theta), \mbox{~for~}0\leq j\leq T.
\end{aligned}
\end{eqnarray}

{\em Step 2:} The induced coefficient map is given by
\begin{equation}\label{exRRC:cm}
\mathcal{C}: \theta \to \Big( f_{T+1}(\theta),\cdots,f_0(\theta),g_{T}(\theta),\cdots,g_0(\theta) \Big).
\end{equation}

{\em Step 3:}  Now it will be shown that the coefficient map \eqref{exRRC:cm} is one-to-one. First it is noted that it can be written as a composition of two functions
\[
\mathcal{C}(\theta) = \mathcal{L} \circ \cR (\theta),
\]
where
\[
\cR:\theta \to \Big(d(\theta), b_1(\theta), a_{1,T}(\theta), \cdots, a_{1,0}(\theta)\Big),
\]
where the components of the vector-function on the right-hand side are given in \eqref{gen-ecm-params}, and
\begin{align*}
\cL:&(d,b_1,a_{1,T},\cdots,a_{1,0}) \to \\ & \hspace{3em}
\Big( f_{T+1},\ldots, f_0,g_{T},\cdots, g_0\Big),
\end{align*}
where the components of this vector-function are given in \eqref{exRRC:tfcoef}.


In order to show that $\mathcal{C}$ is one-to-one it is enough to prove that the maps $\cR$ and $\cL$, defined above, are both one-to-one.

\smallskip

It is claimed that $\cL$ is invertible and thus one-to-one. The inverse map $\cL^{-1}$ is obtained by solving the following set of algebraic equations recursively:
\begin{eqnarray}\label{exRRC:gd}
\begin{aligned}
& d= f_{T+1}\\
& a_{1,0}=-g_{T}\\
& b_1=f_{T}-g_{T}f_{T+1}\\
& a_{1,j}=-g_{T-j}, \mbox{~for~}1\leq j\leq T.
\end{aligned}
\end{eqnarray}

\medskip

\begin{figure}[h]
\centering
\includegraphics[scale=.4]{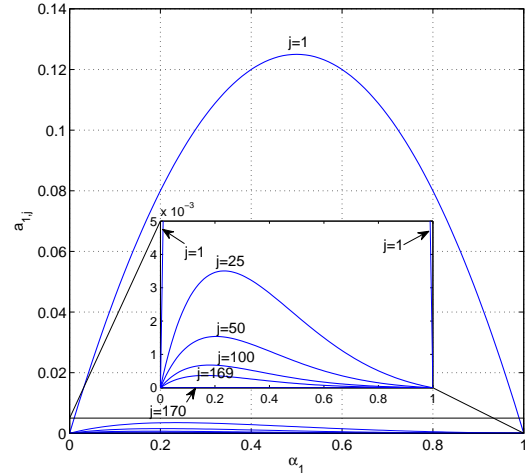}
\caption{$a_{1,j}$ versus $\alpha_1$ within the range $0 \leq \alpha_1 \leq 1$. It is seen that $a_{1,1}$ is symmetric, and that $a_{1,j}\approx 0$ for large $j$. }
\label{Fig:alphaja1}
\end{figure}

Now it is claimed that $\cR$ is one-to-one.  First it is shown that $\alpha_1$ can be uniquely determined by solving the second equation in \eqref{gen-ecm-params} for $i=1$, i.e.
\begin{equation}\label{ex:a11}
a_{1,j}(\alpha_1)= (-1)^{j}\binom{\alpha_1}{j+1}\mbox{~for~} 1 \leq j \leq T,
\end{equation}
where again, since $\alpha_1\in(0,1)$, the generalised binomial coefficient is computed using formula \eqref{Gamma_func}. Note that the function $a_{1,1}(\alpha_1)$ is the only function in  \eqref{ex:a11}, which is  symmetric with respect to the vertical line $\alpha_1=0.5$. In fact it can be shown that for any $j\ne 1$, $a_{1,j}(\alpha_1)$ is nonsymmetric (see Figure~\ref{Fig:alphaja1}). The figure shows that $a_{1,j}\approx 0$ for large $j$.

\begin{figure}[h]
\centering
\includegraphics[scale=.4]{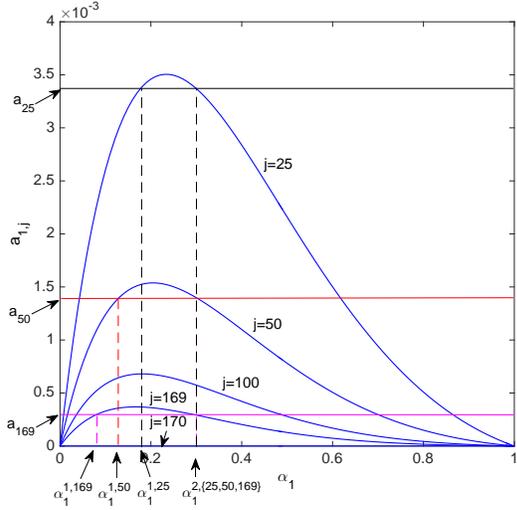}
\caption{A system of at least two equations $a_{1,r}(\alpha_1) = a_{r}$, $a_{1,q}(\alpha_1) = a_{q}$, where $r\ne q$ and 
$1\leq \{r,q\} \leq T$ is sufficient in order to uniquely determine the fractional-order $\alpha_1$.}
\label{Fig:alphaja2}
\end{figure}

Since each of the function \eqref{ex:a11} is unimodal its preimage does not allow for unique determination of the fractional-order $\alpha_1$. However, as shown in Figure \ref{Fig:alphaja2}, the preimage of at least two distinct coefficients $a_{1,r}(\alpha_1)$ and $a_{1,q}(\alpha_1)$, with $1\leq \{r,q\} \leq T$ are sufficient to uniquely determine $\alpha_1$. Using the notation used in Figure \ref{Fig:alphaja2}, each equation $a_{1,j}(\alpha_1) = a_{j}$ gives two solutions  $\alpha_1^{1,j}$ and $\alpha_1^{2,j}$, where $j \in \{1,\cdots,T\}$. Therefore, two distinct coefficients $a_{1,r}$ and $a_{1,q}$ are used to form a system of algebraic equations, whose common solution is the fractional-order $\alpha_1$. This has been shown in Figure \ref{Fig:alphaja2} for $j=25, 50, 169$ with a true value of $\alpha_1=0.3$ (for the clarity of the presentation, $a_{1,1}$ is not plotted in this figure).

Having determined $\alpha_1$ the inverse map $\cR^{-1}$ can be obtained by solving the following set of equations
\[
R_\infty =d, \qquad  C_1=\frac{T_s^{\alpha_1}}{b_1},\qquad R_1=\frac{T_s^{\alpha_1}}{(\alpha_1-a_{1,0})C_1},
\]
which completes the claim and the proof of the proposition.  \hfill{$\square$}

Practical identification of this model using both synthetic and real data has been studied in \cite{Alavi2015}.

%

\medskip

\subsubsection{$R_\infty - R_1||\frac{1}{C_1 s^{\alpha_1}} - R_2||\frac{1}{C_2 s^{\alpha_2}}$ circuit}
\label{exRRCRC}
\begin{figure}[t]
\centering
\includegraphics[scale=.7]{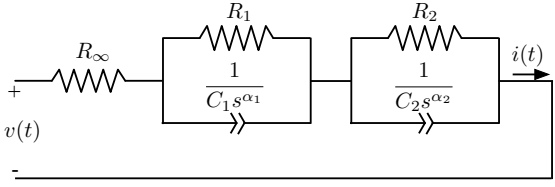}
\caption{The $R_\infty - R_1||\frac{1}{C_1 s^{\alpha_1}} - R_2||\frac{1}{C_2 s^{\alpha_2}}$ circuit.}
\label{Fig:fracR-RC-RC}
\end{figure}

In this part, a more complex model is considered as shown in Figure \ref{Fig:fracR-RC-RC}. The circuit includes two parallel pairs of R-CPE. If $R_2=\infty$, the CPE models the Warburg term.

\medskip

\begin{Proposition}
{\it The FO-ECM Figure~\ref{Fig:fracR-RC-RC} with the parameter vector $\theta=[R_\infty, R_1,  R_2,  C_1,  C_2,  \alpha_1, \alpha_2]$  is structurally identifiable.}
\end{Proposition}


\smallskip

The above proposition will be shown in three steps.

{\em Step 1:} Using \eqref{TF-gecm}, the transfer function of the circuit is
\begin{align}
\label{exRRCRC:tf1}
H(z,\theta)=d(\theta)+\displaystyle \sum_{i=1}^{2}\frac{b_i(\theta)z^{T}}{z^{T+1}- \sum_{j=0}^{T} a_{i,j}(\theta)z^{T-j}},
\end{align}
with the coefficients given in \eqref{gen-ecm-params}. In a rational function form \eqref{exRRCRC:tf1} is expressed as
\begin{align}
\nonumber
&H(z,\theta)=\\
\label{exRRCRC:tf2} &\frac{f_{2T+2}(\theta)z^{2T+2}+f_{2T+1}(\theta)z^{2T+1}+\cdots+ f_{0}(\theta)}{z^{2T+2}+g_{2T+1}(\theta)z^{2T+1}+\cdots+ g_{0}(\theta)}.
\end{align}

For the identifiability analysis, a number of key coefficients are given by:
\begin{eqnarray}\label{exRRCRC:tfcoef}
\begin{aligned}
 f_{2T+2}(\theta)&=d\\
 f_{2T+1}(\theta)&=b_1+b_2+d(-a_{1,0}-a_{2,0})\\
 f_{2T}(\theta)&=-b_1a_{2,0}-b_2a_{1,0}+\\&\hspace{3em}d(a_{1,0}a_{2,0}-a_{1,1}-a_{2,1})\\
 g_{2T+1}(\theta)&=-a_{1,0}-a_{2,0}\\
 g_{2T}(\theta)&=-a_{1,1}-a_{2,1}+a_{1,0}a_{2,0}\\
 g_{2}(\theta)&=a_{1,T} a_{2,T-2} + a_{1,T-1} a_{2,T-1} +\\&\hspace{3em} a_{1,T-2} a_{2,T}\\
 g_{1}(\theta)&=a_{1,T}a_{2,T-1}+a_{1,T-1}a_{2,T}\\
 g_{0}(\theta)&=a_{1,T}a_{2,T}.
\end{aligned}
\end{eqnarray}

{\em Step 2:} The induced coefficient map is given by
\[
\mathcal{C}: \theta \to \Big( f_{2T+2}(\theta),\cdots,f_0(\theta),g_{2T+1}(\theta),\cdots,g_0(\theta) \Big).
\]

{\em Step 3:} Now it will be shown that the above coefficient map is identifiable, i.e. is finitely many-to-one. In order to show the identifiability the following lemma will be used.

\begin{Lemma}\label{ex2:lem}
{\em Consider $g_0$, $g_1$ and $g_2$, the three components of the coefficient map $\mathcal{C}$ expressed in \eqref{exRRCRC:tfcoef}. Then following relations hold:}
\begin{equation}\label{ex2:lemeq}
\begin{aligned}
&g_1 + g_0(T+1)\left(\frac{1}{\alpha_1-T}+\frac{1}{\alpha_2-T}\right) = 0\\
&g_2  - g_0(T+1)(\hat{a}+\hat{b}+\hat{c}) = 0,
\end{aligned}
\end{equation}
where
\begin{align*}
&\hat{a}=\frac{T}{(\alpha_2-T)(\alpha_2-T+1)}\\
&\hat{b}=\frac{(T+1)}{(\alpha_1-T)(\alpha_2-T)}\\
&\hat{c}=\frac{T}{(\alpha_1-T)(\alpha_1-T+1)}.
\end{align*}
\hfill{$\square$}
\end{Lemma}

\noindent{\it Proof:}  To prove the lemma, first it is noted that the relation between $a_{i,j}$ and $a_{i,j+1}$ is the following:

\begin{align*}
a_{i,j+1}&= -(-1)^{j+2}\binom{\alpha_i}{j+2}\\  &=-(-1)^{j+2}\frac{\Gamma(\alpha_i+1)}{\Gamma(j+3)\Gamma(\alpha_i-j-1)}\\ 
&=-1\times -(-1)^{j+1}\frac{\Gamma(\alpha_i+1)}{(j+2)\Gamma(j+2)\frac{\Gamma(\alpha_i-j)}{(\alpha_i-j-1)}}\\&=-\frac{(\alpha_i-j-1)}{j+2}a_{i,j}.
\end{align*}

Therefore, $g_1$ and $g_2$ can be re-written as follows:
\begin{align*}
\nonumber
g_{1} &=a_{1,T}a_{2,T-1}+a_{1,T-1}a_{2,T}=\\
\nonumber & a_{1,T} \frac{-a_{2,T}(T+1)}{\alpha_2-(T-1)-1}+\frac{-a_{1,T}(T+1)}{\alpha_1-(T-1)-1}a_{2,T}\\
& =-g_0(T+1)\left(\frac{1}{\alpha_2-T}+\frac{1}{\alpha_1-T}\right),
\end{align*}

\begin{align*}
\nonumber 
g_{2}&=a_{1,T} a_{2,T-2} + a_{1,T-1} a_{2,T-1} + a_{1,T-2} a_{2,T}\\
\nonumber &=a_{1,T}\frac{-a_{2,T-1}(T-1+1)}{\alpha_2-(T-2)-1}+ \\
\nonumber &\hspace{3em} \frac{-a_{1,T}(T+1)}{\alpha_1-(T-1)-1}\frac{-a_{2,T}(T+1)}{\alpha_2-(T-1)-1}+\\
\nonumber & \hspace{3em}\frac{-a_{1,T-1}(T-1+1)}{\alpha_1-(T-2)-1}a_{2,T}\\
\nonumber & = a_{1,T}\frac{-\frac{-a_{2,T}(T+1)}{\alpha_2-(T-1)-1}(T-1+1)}{\alpha_2-(T-2)-1} +\\ \nonumber & \hspace{3em} \frac{-a_{1,T}(T+1)}{\alpha_1-(T-1)-1}\frac{-a_{2,T}(T+1)}{\alpha_2-(T-1)-1} + \\
\nonumber  &\hspace{5em}\frac{-\frac{-a_{1,T}(T+1)}{\alpha_1-(T-1)-1}(T-1+1)}{\alpha_1-(T-2)-1}a_{2,T}\\
&=g_0(T+1)(\hat{a}+\hat{b}+\hat{c}).
\end{align*} \hfill{$\square$}

By the result of Lemma~\ref{ex2:lem} the exponent values $\alpha_1$ and $\alpha_2$ are determined by solving  \eqref{ex2:lemeq}. It can be easily verified (using any computer algebra package, e.g.\ Mathematica) that these algebraic equations admit only two real solutions, which are permuted with respect to each other. Since the mathematical expressions of the solutions are quite cumbersome, they are not provided explicitly in this paper. 

Finally, it is noted that from the first equation in \eqref{exRRCRC:tfcoef}, the ohmic resistor $R_\infty$ is estimated trough the coefficient $f_{2T+2}$, i.e., $R_\infty=f_{2T+2}$. Then $a_{i,j}$'s are calculated recursively, backward from $j=T$ to $j=1$ for $i=1,2$. By obtaining $a_{1,1}$ and $a_{2,1}$, the parameters $a_{1,0}$ and $a_{2,0}$ can be computed from the coefficients $g_{2T}$ and $g_{2T+1}$ in \eqref{exRRCRC:tfcoef}. Then $b_1$ and $b_2$ are computed from the coefficients  $f_{2T}$ and $f_{2T+1}$ in \eqref{exRRCRC:tfcoef}. And, for a given solution of \eqref{ex2:lemeq}  (there are two real solutions) parameters $C_1$ and $C_2$ are then obtained by using $b_i=T_s^{\alpha_i}/C_i$, and finally $R_1$ and $R_2$ are calculated from $a_{1,0}$ and $a_{2,0}$, respectively.  This completes the Step~3 and the proof of the proposition. \hfill{$\square$}

\section{Conclusions}
A method was proposed for the structural identifiability analysis of fractional order (FO) systems based on the concept of  coefficient map. The method is applicable to both commensurate and non-commensurate models and was applied to determine the structural identifiability of battery fractional-order equivalent circuit models (FO-ECMs). This study has shown that a battery FO-ECM is structurally  identifiable, and the global identifiability was proved for the FO-ECM with a single  constant phase element.

\section*{Acknowledgements}
Authors would like to acknowledge fundings from the UK's Engineering and Physical Sciences Research Council (EPSRC), S.M.M. Alavi and D.A. Howey under grant EP/K503769/1, A. Mahdi and S.J. Payne under grant EP/K036157/1 and P.E. Jacob under grant EP/K009362/1.
          
\bibliographystyle{IEEEtran}
\bibliography{Part1_frac_ident_final2.bbl}
\end{document}